\font\largebf=cmbx10 scaled\magstep2
\font\tenmsa=msam10
\font\tenmsb=msbm10

\def\all{\hbox{for all}}
\def\and{\hbox{and}}

\def\bra#1#2{\langle#1,#2\rangle}
\def\Bra#1#2{\big\langle#1,#2\big\rangle}

\def\cite#1\endcite{[#1]}
\def\Cite#1\endcite{\big[#1\big]}

\def\dom{\hbox{\rm dom}}

\def\f#1#2{{#1 \over #2}}

\def\half{\ts\f12}

\def\lr{\Longrightarrow}

\def\NQ{{\cal N}_q}

\def\on{\hbox{on}}
\def\PC{{\cal PC}}
\def\PCLSC{{\cal PCLSC}}

\def\phi{\varphi}

\def\qed{\hfill\hbox{\tenmsa\char03}}
\def\qlr{\quad\lr\quad}

\def\r{\hbox{\rm I \hskip - .5em R}}
\def\r{\hbox{\tenmsb R}}
\def\rbar{\r\cup\{\infty\}}
\def\rbar{\,]{-}\infty,\infty]}

\long\def\slant#1\endslant{{\sl#1}}

\def\supn{\sup\nolimits}

\def\ts{\textstyle}

\def\xbra#1#2{\lfloor#1,#2\rfloor}

\def\Xbra#1#2{\big\lfloor#1,#2\big\rfloor}

\def\defSection#1{}
\def\defCorollary#1{}
\def\defDefinition#1{}
\def\defExample#1{}
\def\defLemma#1{}
\def\defNotation#1{}
\def\defProblem#1{}
\def\defRemark#1{}
\def\defTheorem#1{}
\def\locno#1{}
\def\meqno#1{\eqno(#1)}
\def\nmbr#1{}
\def\Proof{\medbreak\noindent{\bf Proof.}\enspace}
\def\Proo{{\bf Proof.}\enspace}
\def\Signoff{}
\def \INTsec{1}
\def \SSDsec{2}
\def \SSDBdef{2.1}
\def \SSDBone{1}
\def \SSDBtwo{2}
\def \SSDBthree{3}
\def \SSDBfour{4}
\def \SSDBex{2.2}
\def \EEex{2.3}
\def \TRANSsec{3}
\def \AT{5}
\def \TRANSlem{3.1}
\def \TRANSone{6}
\def \TRANStwo{7}
\def \RTRlem{3.2}
\def \DOMthm{3.3}
\def \DOMone{8}
\def \DOMtwo{9}
\def \DOMrem{3.4}
\def \LINsec{4}
\def \SSDlem{4.1}
\def \Hlem{4.2}
\def \Hone{10}
\def \PERPlem{4.3}
\def \PERPthm{4.4}
\def \BBcor{4.5}
\def \BB{1}
\def \ML{2}
\def \MLP{3}
\def \RTRFENCHEL{4}
\def \HBM{5}
\def \SSDMON{6}
\def \SSDB{7}
\def \YAO{8}
\magnification 1200
\headline{\ifnum\folio=1
{\hfil{\largebf A Brezis--Browder theorem for SSDB spaces}\hfil}
\else\centerline{\rm {\bf A Brezis--Browder theorem for SSDB spaces}}\fi}
\medskip
\centerline{by}
\medskip
\centerline{Stephen Simons}
\medskip
\centerline{\bf Abstract}
\medskip\noindent
In this paper, we show how the Brezis-Browder theorem for maximally monotone multifunctions with a linear graph on a reflexive Banach space, and a consequence of it due to Yao, can be generalized to SSDB spaces.
\defSection \INTsec
\bigbreak
\centerline{\bf \INTsec\quad Introduction}
\medskip
\noindent This note is about generalizations of the following results of Brezis-Browder and Yao: \slant Let $E$ be a nonzero reflexive Banach space with topological dual $E^*$ and $A$ be a norm--closed monotone linear subspace of $E \times E^*$.   Then $A$ is maximally monotone $\iff A^*$ is monotone $\iff A^*$ is maximally monotone.\endslant\  In Theorem \PERPthm, we show that these results can be successfully generalized to SSDB space (see Section \SSDsec).   The original results of Brezis--Browder and Yao appear in Corollary \BBcor.   Theorem \PERPthm\ depends on a new transversality result, which appears in Theorem \DOMthm.  
\defSection \SSDsec
\medbreak
\centerline{\bf \SSDsec\quad SSDB spaces}
\medskip
\noindent All vector spaces in the paper will be \slant real\endslant.   We will use the standard notation of convex analysis for ``Fenchel conjugate'' and ``subdifferential'' in a Banach space.
\defDefinition \SSDBdef
\medbreak
\noindent
{\bf Definition \SSDBdef.}\enspace We will say that $\big(B,\xbra\cdot\cdot,q,\|\cdot\|,\iota\big)$ is a \slant symmetrically self--dual Banach space (SSDB space)\endslant\ if $B$ is a nonzero real vector space, $\xbra\cdot\cdot\colon B \times B \to \r$ is a symmetric bilinear form,\quad for all $b \in B$, $q(b) := \half\xbra{b}{b}$\quad   (``$q$'' stands for ``quadratic''), $\big(B,\|\cdot\|\big)$ is a Banach space, and $\iota$ is a linear isometry from $B$ onto $B^*$ such that, for all $b,c \in B$, $\Bra{b}{\iota(c)} = \xbra{b}{c}$.   It is easily seen that if $\big(B,\xbra\cdot\cdot,q,\|\cdot\|,\iota\big)$ is a SSDB space then so also is $\big(B,-\xbra\cdot\cdot,-q,\|\cdot\|, -\iota\big)$.   Clearly, for all $b,c \in B$,
$$\left.\eqalign{\big|q(b) - q(c)\big| &= \half\big|\xbra{b}{b} - \xbra{c}{c}\big| = \half\big|\xbra{b - c}{b + c}\big| = \half\big|\Bra{b - c}{\iota(b + c)}\big|\cr
&\le \half\|b - c\|\big\|\iota(b + c)\big\| =  \half\|b - c\|\|b + c\|.}\right\}\meqno\SSDBone$$
We refer the reader to Examples \SSDBex\ and \EEex\ below and \cite\MLP\endcite\ for various examples.   There is also a construction that can be used for producing more pathological examples in \cite\SSDMON, Remark 6.7, p.\ 20\endcite.
\smallbreak
Let $\big(B,\xbra\cdot\cdot,q,\|\cdot\|,\iota\big)$ be an SSDB space and $A \subset B$.   We say that $A$ is \slant$q$--positive\endslant\ if $A \ne \emptyset$ and\quad $b,c \in A \lr q(b - c) \ge 0$.\quad We say that $A$ is \slant$q$--negative\endslant\ if\quad $A \ne \emptyset$ and\quad $b,c \in A \lr q(b - c) \le 0$.   $q$--negativity is equivalent to $(-q)$--positivity.   We say that $A$ is \slant maximally $q$--positive\endslant\ if $A$ is $q$--positive and $A$ is not properly contained in any other $q$--positive set.   In this case,
$$b \in B \lr \inf q(A - b) \le 0.\meqno\SSDBtwo$$
Similarly, we say that $A$ is \slant maximally $q$--negative\endslant\ if $A$ is $q$--negative and $A$ is not properly contained in any other $q$--negative set.   Maximal $q$--negativity is equivalent to maximal $(-q)$--positivity.
\medskip
Let $g_0 := \half\|\cdot\|^2$ on $B$.   Then, for all $c \in B$,\quad
${g_0}^*\big(\iota(c)\big) = \half\|\iota(c)\|^2 = \half\|c\|^2 = g_0(c)$.\quad
Consequently, for all $b,c \in B$,
$$\iota(c) \in \partial g_0(b) \iff g_0(b) + {g_0}^*\big(\iota(c)\big) = \Bra{b}{\iota(c)} \iff g_0(b) + g_0(c) = \xbra{b}{c}.
\meqno\SSDBthree$$
We note from (\SSDBone) with $c = 0$ that
$$g_0 + q \ge 0\ \on\ B.\meqno\SSDBfour$$\par
\defExample \SSDBex 
\medbreak
\noindent¨
{\bf Examples \SSDBex.}
\smallbreak
(a)\enspace  If $B$ is a Hilbert space with inner product $(b,c) \mapsto \bra{b}{c}$ and the Hilbert space norm then $\big(B,\xbra\cdot\cdot,q,\|\cdot\|,\iota\big)$ is a SSDB space with $\xbra{b}{c} := \bra{b}{c}$,  $q = g_0$ and $\iota(c) := c$, and every nonempty subset of $B$ is $q$--positive.
\smallbreak
(b)\enspace  If $B$ is a Hilbert space with inner product $(b,c) \mapsto \bra{b}{c}$ and the Hilbert space norm then $\big(B,\xbra\cdot\cdot,q,\|\cdot\|,\iota\big)$ is a SSDB space with $\xbra{b}{c} := -\bra{b}{c}$, $q = -g_0$ and $\iota(c) := -c$, and the $q$--positive sets are the singletons.
\smallbreak
(c)\enspace If $B = \r^3$ under the Euclidean norm and 
$$\Xbra{(b_1,b_2,b_3)}{(c_1,c_2,c_3)} := b_1c_2 + b_2c_1 + b_3c_3,$$
then $\big(B,\xbra\cdot\cdot,q,\|\cdot\|,\iota\big)$ is a SSDB space,\quad$q(b_1,b_2,b_3) = b_1b_2 + \half b_3^2$ and $\iota(c_1,c_2,c_3) := (c_2,c_1,c_3)$.   Here, if $M$ is any nonempty monotone subset of $\r \times \r$ (in the obvious sense) then $M \times \r$ is a $q$--positive subset of $B$.   The set $\r(1,-1,2)$ is a $q$--positive subset of $B$ which is not contained in a set $M \times \r$ for any monotone subset of $\r \times \r$.   The helix $\big\{(\cos\theta,\sin\theta,\theta)\colon \theta \in \r\big\}$ is a $q$--positive subset of $B$, but if
$0 < \lambda < 1$ then the helix $\big\{(\cos\theta,\sin\theta,\lambda\theta)\colon \theta \in \r\big\}$ is not.
%
\defExample \EEex
\medbreak
\noindent
{\bf Example \EEex.}\enspace Let $E$ be a nonzero reflexive Banach space with topological dual $E^*$.   For all $(x,x^*)$ and $(y,y^*) \in E \times E^*$, let\quad $\Xbra{(x,x^*)}{(y,y^*)} := \bra{x}{y^*} + \bra{y}{x^*}$.\quad We norm $E \times E^*$ by $\big\|(x,x^*)\big\| := \sqrt{\|x\|^2 + \|x^*\|^2}$.   Then $\big(E \times E^*,\|\cdot\|\big)^* = (E^* \times E,\|\cdot\|\big)$,\break under the duality\quad $\Bra{(x,x^*)}{(y^*,y)} := \bra{x}{y^*} + \bra{y}{x^*}$.\quad  It is clear from these definitions that\quad $\Bra{(x,x^*)}{(y^*,y)} = \Xbra{(x,x^*)}{(y,y^*)}$.\quad Thus $\big(E \times E^*,\xbra\cdot\cdot,q,\|\cdot\|,\iota\big)$ is a SSDB space with\quad $q(x,x^*) = \half\big[\bra{x}{x^*} + \bra{x}{x^*}\big] = \bra{x}{x^*}$\quad and\quad $\iota(y,y^*) = (y^*,y)$.\quad
\par
We now note that if\quad $(x,x^*), (y,y^*) \in B$\quad then\quad$\bra{x - y}{x^* - y^*} = q(x - y,x^* - y^*) = q\big((x,x^*) - (y,y^*)\big)$.\quad  Thus if $A \subset B$ then $A$ is $q$--positive exactly when $A$ is a nonempty monotone subset of $B$ in the usual sense, and $A$ is maximally $q$--positive exactly when $A$ is a maximally monotone subset of $B$ in the usual sense.
\par
We define the reflection map $\rho_1\colon\ E \times E^* \to E \times E^*$ by $\rho_1(x,x^*) := (-x,x^*)$.   Since $q \circ \rho_1 = -q$, a subset $A$ of $B$ is $q$--positive (resp. maximally $q$--positive) if, and only if, $\rho_1(A)$ is $q$--negative (resp. maximally $q$--negative).  $\rho_1$ and its companion map $\rho_2$ are used in the discussion of an abstract Hammerstein theorem in \cite\HBM, Section 30, pp.\ 123-- 125\endcite\ and \cite\SSDB, Section 7, pp. 13--15\endcite. 
\defSection \TRANSsec
\bigbreak
\centerline{\bf \TRANSsec\quad  A transversality result}
\medskip
\noindent
Let $B = \big(B,\xbra\cdot\cdot,q,\|\cdot\|,\iota\big)$ be a SSDB space.   We write $\PC(B)$ for the set of all proper convex functions from $B$ into $\rbar$ and $\dom\,f$ for the set $\big\{x \in B\colon\ f(x) \in \r\big\}$.    If $f \in \PC(B)$, we write $f^@$ for the Fenchel conjugate of $f$ with respect to the pairing $\xbra\cdot\cdot$, that is to say, for all $c \in B$,\quad$f^@(c) := \supn_B\big[\xbra{\cdot}{c} - f\big]$.   We note then that, for all $c \in B$,
$$f^@(c) = \supn_B\big[\bra{\cdot}{\iota(c)} - f\big] = f^*\big(\iota(c)\big).\meqno\AT$$
We write $\PCLSC(B)$ for the set $\{f \in \PC(B)\colon\ f\ \hbox{is lower semicontinuous on}\ B\}$.
\defLemma \TRANSlem
\medbreak
\noindent
{\bf Lemma \TRANSlem.}\enspace\slant Let $\big(B,\xbra\cdot\cdot\big)$ be a SSDB space, $f \in \PC(B)$ and $c \in B$.   We define\quad $f_c\in \PC(B)$\quad by\quad $f_c := f(\cdot + c) - \xbra{\cdot}{c} - q(c)$.\quad  Then
$$(f_c)^@ = (f^@)_c\meqno\TRANSone$$
and, writing $f_c^@$ for the common value of these two functions,
$$\all\ b,d \in B,\ f_c(b) + f_c^@(d) - \xbra{b}{d} = f(b+ c) + f^@(d + c) - \xbra{b + c}{d + c}.\meqno\TRANStwo$$\endslant
\Proof For all $b \in B$,
$$\leqalignno{{f_c}^@(b)
&= \supn_{d \in B}\big[\xbra{d}{b} + \xbra{d}{c} + q(c) - f(d + c)\big]\cr
&= \supn_{e \in B}\big[\xbra{e - c}{b + c} + q(c) - f(e)\big]\cr
&= \supn_{e \in B}\big[\xbra{e}{b + c} - \xbra{c}{b} - f(e)\big] - q(c)\cr
&= f^@(b + c) - \xbra{c}{b} - q(c) = (f^@)_c(b),}$$
which gives (\TRANSone), and (\TRANStwo) follows since
$$\eqalignno{f_c(b) + f_c^@(d) - \xbra{b}{d} &= f(b + c) - \xbra{b}{c} - q(c) + f^@(d + c) - \xbra{d}{c} - q(c) - \xbra{b}{d}\cr
&= f(b + c) + f^@(d + c) - \xbra{b}{c} - \xbra{c}{c} - \xbra{d}{c} - \xbra{b}{d}\cr
&= f(b+ c) + f^@(d + c) - \xbra{b + c}{d + c}.&\qed}$$\par
Theorem \DOMthm\ is the main result of this section, and uses the following consequence of Rockafellar's formula for the subdifferential of a sum \big(see \cite\RTRFENCHEL, Theorem 3, p.\ 86\endcite\big).  
\defLemma \RTRlem
\medbreak
\noindent
{\bf Lemma \RTRlem.}\enspace\slant Let $X$ be a nonzero normed
space, $f\colon\ X \mapsto \rbar$ be convex on $X$, finite at a point of $X$, and $g\colon\ X \mapsto \r$ be convex and continuous.   Then $\partial (f + g) = \partial f + \partial g$.\endslant
\defTheorem \DOMthm
\medbreak
\noindent
{\bf Theorem \DOMthm.}\enspace\slant Let $B = \big(B,\xbra\cdot\cdot,q,\|\cdot\|,\iota\big)$ be a SSDB space, $f \in \PCLSC(B)$ and, whenever $b,d \in B$,
$$f(b) + f^@(d) = \xbra{b}{d} \qlr q(b) + q(d) \le \xbra{b}{d}.\meqno\DOMone$$
Let\quad $\NQ(g_0) := \big\{b \in B\colon\ g_0(b) + q(b) = 0\big\}$.\quad Then\quad $\dom\,f - \NQ(g_0) = B$.\endslant
\Proof Let $c$ be an arbitrary element of $B$.   Now\quad $f_c + g_0$\quad is convex, and \quad $\big(f_c + g_0\big)(b) \to \infty$ as $\|b\| \to \infty$,\quad i.e., $f_c + g_0$ is ``coercive''.   Since $g_0$ is continuous and convex, $f_c + g_0 \in \PCLSC(B)$, hence $f_c + g_0$ is $w(B,B^*)$--lower semicontinuous.   From \cite\ML, Proposition 1.3, p.\ 892\endcite, $B$ is reflexive, and so $f_c + g_0$ attains a minimum at some $b \in B$. Thus $\partial\big(f_c + g_0\big)(b) \ni 0$.   Since $g_0$ is continuous, Lemma \RTRlem\ implies that $\partial f_c(b) + \partial g_0(b) \ni 0$, and so there exists $d^* \in \partial f_c(b)$ such that $-d^* \in \partial g_0(b)$.   Since $\iota$ is surjective, there exists $d \in B$ such that $\iota(d) = d^*$.   Then, from (\SSDBthree),
$$g_0(b) + g_0(d) = -\xbra{b}{d}\meqno\DOMtwo$$
and, using (\AT),\quad $f_c(b) + f_c^@(d) = f_c(b) + {f_c}^*(d^*) = \bra{b}{d^*} = \xbra{b}{d}$.\quad   Combining this with (\TRANStwo),\quad $f(b + c) + f^@(d + c) = \xbra{b + c}{d + c}$.\quad   Since this gives\quad $f(b + c) \in \r$,\quad $c \in \dom\,f - b$.\quad   From (\DOMone),\quad  $q(b + c) + q(d +c ) \le \xbra{b + c}{d + c}$,\quad or equivalently,\quad $q(b) + q(d) \le \xbra{b}{d}$.\quad   Adding this to (\DOMtwo),\quad $(g_0 + q)(b) + (g_0 + q)(d) \le 0$.\quad   From (\SSDBfour),\quad $(g_0 + q)(b) = 0$,\quad that is to say,\quad $b \in \NQ(g_0)$.\quad   Thus\quad $c \in \dom\,f - \NQ(g_0)$,\quad as required.\qed
\defRemark \DOMrem
\medbreak
\noindent
{\bf Remark \DOMrem.}\enspace While Theorem \DOMthm\ is adequate for the application in Lemma \PERPlem, if we think of $\partial f$ as a multifunction from $B$ into $B^*$, then the proof above actually establishes that \quad $D(\partial f) - \NQ(g_0) = B$.
\defSection \LINsec
\medbreak
\centerline{\bf \LINsec\quad  Linear $q$--positive sets}
\medskip
\noindent
We now come to the main topic of this paper:  linear $q$--positive sets.   Let $\big(B,\xbra\cdot\cdot,q,\|\cdot\|,\iota\big)$ be an SSDB space and $A$ be a \slant linear subspace\endslant\ of $B$.   Then we write $A^0$ for the linear subspace  $\big\{b \in B\colon\ \xbra{A}{b} = \{0\}\big\}$ of $B$.  
\defLemma \SSDlem
\medbreak
\noindent
{\bf Lemma \SSDlem.}\enspace\slant Let $B = \big(B,\xbra\cdot\cdot,q,\|\cdot\|,\iota\big)$ be a SSDB space and $A$ be a maximally $q$--positive linear subspace of $B$.   Then $A^0$ is a $q$--negative linear subspace.\endslant
\Proof If $p \in A^0$ then\quad $\inf q(A - p) = \inf q(A) + q(p) = q(p)$,\quad and so (\SSDBtwo) gives $q(p) \le 0$.   If $b,c \in A^0$ then $b - c \in A^0$ and so $q(b - c) \le 0$.\qed
%
\defLemma \Hlem
\medskip
\noindent
{\bf Lemma \Hlem.}\enspace\slant Let $B = \big(B,\xbra\cdot\cdot,q,\|\cdot\|,\iota\big)$ be a SSDB space and $A$ be a norm--closed $q$--positive linear subspace of $B$.   Define the function\quad $q_A\colon\ B \to \,\rbar$\quad by\quad $q_A := q$ on $A$ and $q_A := \infty$ on $B \setminus A$.   Then $q_A \in \PCLSC(B)$ and
$$q_A(b) + {q_A}^@(d) = \xbra{b}{d} \lr b - d \in A^0.\meqno\Hone$$\endslant
\Proo Suppose that $a,c \in A$ and $\lambda \in \,]0,1[\,$.   Then
$$\lambda q(a) + (1 - \lambda)q(c) - q\big(\lambda a + (1 - \lambda)c\big) = \lambda(1 - \lambda)q(a - c) \ge 0.$$
It is easily seen that this implies the convexity of $q_A$.  \big(See \cite\HBM, Lemma 19.7, pp.\ 80--81\endcite.\big)  We know from (\SSDBone) that $q$ is continuous.  Since $A$ is norm--closed in $B$, $q_A \in \PCLSC(B)$.
\smallbreak
We now establish (\Hone).   Let\quad $q_A(b) + {q_A}^@(d) = \xbra{b}{d}$, $a \in A$ and $\lambda \in \r$.\quad  Then clearly\quad $b \in A$ and $b + \lambda a \in A$,\quad and so $\xbra{b + \lambda a}{d} - \xbra{\lambda a}{b} - \lambda^2q(a) = q(b) + \xbra{b + \lambda a}{d} - q(b + \lambda a) = q_A(b) + \xbra{b + \lambda a}{d} - q_A(b + \lambda a)\le q_A(b) + {q_A}^@(d) = \xbra{b}{d}$.   Thus $\lambda^2q(a) + \lambda\xbra{b - d}{a} \ge 0$.   Since this inequality holds for all $\lambda \in \r$, $\xbra{b - d}{a} = 0$.   Since this equality holds for all $a \in A$, we obtain (\Hone).\qed
\defLemma \PERPlem
\medbreak
\noindent
{\bf Lemma \PERPlem.}\enspace\slant Let $B = \big(B,\xbra\cdot\cdot,q,\|\cdot\|,\iota\big)$ be a SSDB space, $A$ be a norm--closed $q$--positive linear subspace of $B$, and $A^0$ be $q$--negative.   Then:
\par\noindent
{\rm(a)}\enspace $A - \NQ(g_0) = B$.
\par\noindent
{\rm(b)}\enspace $A$ is maximally $q$--positive.\endslant
\Proof(a)\enspace It is clear from (\Hone) that $q_A(b) + {q_A}^@(d) = \xbra{b}{d} \lr q(b - d) \le 0$, that is to say, $q_A(b) + {q_A}^@(d) = \xbra{b}{d} \lr q(b) + q(d) \le \xbra{b}{d}$.  (a) now follows easily from Theorem \DOMthm\ with $f = q_A$.
\par
(b)\enspace Now suppose that $c \in B$ and\quad $A \cup \{c\}$\quad is $q$--positive.   From (a), there exists $a \in A$ such that\quad $a - c \in \NQ(g_0)$.\quad   Thus\quad $\half\|a - c\|^2 + q(a - c) = 0$. \quad   Since\quad $A \cup \{c\}$\quad is $q$--positive,\quad $q(a - c) \ge 0$,\quad and so\quad $\half\|a - c\|^2 \le 0$,\quad from which\quad $c = a \in A$.\quad   This completes the proof of (b).\qed
\medbreak
Our next result is suggested by Yao, \cite\YAO, Theorem 2.4, p.\ 3\endcite.
\defTheorem \PERPthm
\medbreak
\noindent
{\bf Theorem \PERPthm.}\enspace\slant  Let $B = \big(B,\xbra\cdot\cdot,q,\|\cdot\|,\iota\big)$ be a SSDB space and $A$ be a norm--closed $q$--positive linear subspace \big(of $\big(B,\xbra\cdot\cdot,q,\|\cdot\|,\iota\big)$\big).   Then:
\smallbreak\noindent
{\rm(a)}\enspace $A$ is maximally $q$--positive if, and only if, $A^0$ is $q$--negative.
\smallbreak\noindent
{\rm(b)}\enspace $A$ is maximally $q$--positive if, and only if, $A^0$ is maximally $q$--negative.\endslant
\Proof (a) is clear from Lemmas \SSDlem\ and \PERPlem(b).   ``If'' in (b) is immediate from ``If'' in (a).   Conversely, let us suppose that $A$ is maximally $q$--positive.   We know already from (a) that $A^0$ is $q$--negative, and it only remains to prove the maximality.   Since $A$ is norm--closed, it follows from standard functional analysis and the surjectivity of $\iota$ that $A = (A^0)^0$.  Now $A^0$ is $(-q)$--positive and norm--closed.   Furthermore, $(A^0)^0 = A$ is $q$--positive and thus $(-q)$--negative.   If we now apply Lemma \PERPlem(b), with $A$ replaced by $A^0$ and $q$ replaced by $-q$, we see that $A^0$ is maximally $(-q)$--positive, that is to say, maximally $q$--negative, which completes the proof of (b).\qed          
\medskip
Now let $A$ be a monotone linear subspace of $E \times E^*$.   Then the linear subspace $A^*$ of $E \times E^*$ is defined by:\quad  $(x,x^*) \in A^* \iff \hbox{for all}\ (a,a^*) \in A,\ \bra{x}{a^*} = \bra{a}{x^*}$.\quad It is clear then that $A^0 = \rho_1\big(A^*\big)$.   
Corollary \BBcor(a) below appears in Brezis--Browder \cite\BB, Theorem 2, pp.\ 32--33\endcite, and Corollary \BBcor(b) in Yao, \cite\YAO, Theorem 2.4, p.\ 3\endcite. 
\defCorollary \BBcor
\medbreak
\noindent
{\bf Corollary \BBcor.}\enspace\slant Let $E$ be a nonzero reflexive Banach space with topological dual $E^*$ and $A$ be a norm--closed monotone linear subspace of $E \times E^*$.   Then:
\smallbreak\noindent
{\rm(a)}\enspace $A$ is maximally monotone if, and only if, $A^*$ is monotone
\smallbreak
\noindent
{\rm(b)}\enspace $A$ is maximally monotone if, and only if, $A^*$ is maximally monotone.\endslant
\Proof These results follow from Theorem \PERPthm\ and the comments in Example \EEex.\qed
\bigbreak
\centerline{\bf References}
\medskip
\nmbr\BB
\item{[\BB]} H. Brezis and F. E. Browder, {\sl  Linear maximal monotone operators and singular nonlinear integral equations of Hammerstein,  type}, Nonlinear analysis (collection of papers in honor of Erich H. Rothe),  pp. 31--42. Academic Press, New York, 1978.
\nmbr\ML
\item{[\ML]} J-E. Mart\'\i nez-Legaz, \slant On maximally $q$--positive sets\endslant,  J. of Convex Anal., {\bf 16} (2009), 891--898.
\nmbr\MLP
\item{[\MLP]} J-E. Mart\'\i nez-Legaz, \slant Private communication\endslant.
\nmbr\RTRFENCHEL
\item{[\RTRFENCHEL]} R. T. Rockafellar, \slant Extension of Fenchel's
duality theorem for convex functions\endslant, Duke Math. J. {\bf33}
(1966),  81--89.
\nmbr\HBM
\item{[\HBM]}S. Simons, \slant From Hahn--Banach to monotonicity\endslant,  Lecture Notes in Mathematics, {\bf 1693},\break second edition, (2008), Springer--Verlag.
\nmbr\SSDMON
\item{[\SSDMON]}-----, \slant Banach SSD spaces and classes of monotone sets\endslant, http://arxiv/org/abs/\break 0908.0383v3, posted July 5, 2010.
\nmbr\SSDB
\item{[\SSDB]}-----, \slant SSDB spaces and maximal monotonicity\endslant, http://arxiv/org/abs/1001.0064v2,\break posted June 17, 2010.
\nmbr\YAO 
\item{[\YAO]} L. Yao, \slant The Brezis--Browder theorem revisited and properties of Fitzpatrick functions of order $n$\endslant, http://arxiv.org/abs/0905.4056v1, posted May 25, 2009.
\Signoff
\bigskip
Department of Mathematics\par
University of California\par
Santa Barbara\par
CA 93106-3080\par
U. S. A.\par
email:  simons@math.ucsb.edu
\bye